%
%

\magnification=1200


\font\AAA=cmr14 at 12pt
\font\BBB=cmr14 at 8pt

\overfullrule=0in

\def\boxit#1{\hbox{\vrule
 \vtop{%
  \vbox{\hrule\kern 2pt %
     \hbox{\kern 2pt #1\kern 2pt}}%
   \kern 2pt \hrule }%
  \vrule}}

\def\cn{\bbc^n}
\def\ss{\subset}

\def\dist{{\rm dist}}

\def\log{{\rm log}}

\def\supp{{\rm supp}}

\def\fa{\qquad{\rm for \ all\ \ }}


\def\Theorem#1{\medskip\noindent {\AAA T\BBB HEOREM \rm #1.}}

\def\Cor#1{\medskip\noindent {\AAA C\BBB OROLLARY \rm #1.}}
\def\Lemma#1{\medskip\noindent {\AAA L\BBB EMMA \rm  #1.}}

\def\Note#1{\medskip\noindent {\AAA N\BBB OTE \rm  #1.}}

\def\pf{\medskip\noindent {\bf Proof.}\ }
\def\qed{\hfill  $\vrule width5pt height5pt depth0pt$}

\def\cv{{\cal V}}   \def\cp{{\cal P}}   
\def\cw{{\cal W}}   
\def\ce{{\cal E}}   
   
\def\cs{{\cal S}}   

\def\cp{{\cal P}}

\def\vf{\varphi}

\def\wh{\widehat}

\def\and{\qquad {\rm and} \qquad}

\def\bbc{{\bf C}}

\def\bbp{{\bf P}}

\def\a{\alpha}

\def\d{\delta}

\def\g{\gamma}

\def\l{\lambda}

\def\z{\zeta}

\def\L{\Lambda}
\def\G{\Gamma}
\def\O{\Omega}

\font\titfont=cmr10 at 12 pt

\def\z{\zeta}

\centerline{\titfont  ON THE COMPLEMENT OF THE  PROJECTIVE HULL  IN C$^n$.  }
 
\vskip .2in
\centerline{\titfont  Blaine Lawson$^*$ and John Wermer}
\vglue .9cm
\smallbreak\footnote{}{ $ {} \sp{ *}{\rm Partially}$  supported by
the N.S.F. }

\centerline{\bf ABSTRACT} \medskip
  \font\abstractfont=cmr10 at 10 pt

{{\parindent= .5in
\narrower\abstractfont \noindent
We prove that if $K$ is a compact subset of an affine variety $\Omega = \bbp^n-D$ (where  $D$ is a projective hypersurface) and if  $K$ is contained in  a closed analytic subvariety  $V \ss \Omega$, then the projective hull $\wh K$  has the property that   $\wh K\cap \O\ss V$.
If $V$ is smooth and  1-dimensional, then $\wh K\cap\O$ is also closed in $\O$.  The result has 
applications to graphs in $\bbc^2$ of functions in the disk algebra.

}}

\vskip .3in

Let $X$ be a compact set in $\cn$.  In [1],  R. Harvey and the  first author introduce a generalization of the  polynomial hull of $X$, which they call the {\sl projective hull of $X$ in $\cn$}, and denote by
$\wh X$.  For $d=1,2,...$ let $\cp_d$  denote the space of all 
polynomials on $\cn$ of degree $\leq d$.  
Then by definition a point $z$ lies in $\wh X$ if there exists a constant $C\geq 1$ such that for all $d$ and for all $P\in\cp_d$ with $\sup_X|P|\leq 1$, we have
$$
|P(z)|\ \leq\ C^d.
$$

Fix a point $a=(a_1,...,a_n)\in\cn - X$.
For $d=1,2,...$ set 
$$
\l_d(a)\ \equiv\ \sup\{ |P(a)| : P\in \cp_d {\ \rm with\ }\|P\|_X\leq 1\}
$$
Note that $a$ lies in the complement of $\wh X$ if and only if for each $M\geq 1$ there exist arbitrarily large $d$ with $\l_d(a)\geq M^d$.

We call a set  $X${\sl non-algebraic} if every polynomial vanishing on $X$ is identically zero.
We restrict our attention to such non-algebraic sets.  We study the function $\l_d$ and, in Theorem 2
below, give a result on the complement of $\wh X$ in certain cases, generalizing Theorem 9.2 in [1].

\Theorem {1}  {\sl  Let $X$ be as above.  
 Fix $M\geq 1$.  Fix $d$ and let $\l_d=\l_d(0)$.  Then $\l_d \geq M^d$ if and only if there exist polynomials $A_1,...,A_n $ on $\cn$ such that $\z_1A_1(\z)+\cdots+\z_nA_n(\z) \in \cp_d$ 
with $|1-\sum_k \z_k A_k|\leq {1\over M^d}$ on $X$.}
\medskip

We use the following Banach distance formula.

\Lemma {1}  {\sl  Let $W$ be  a normed linear space, dim$(\cw)=N$, and let  $\cv$ be a 
subspace of $\cw$ with dim$(\cv)=N-1$.  Fix $x\in \cw-\cv$.  Let $L$ be the linear functional
on $\cw$ with $L=0$ on $\cv$ and $L(x)=1$.  Let $\d = \dist(x, \cv)$  Then}
$$
\|L\|\ =\ {1\over \d}.
$$

\Note{} See the Appendix for a proof. 
\medskip
\noindent
{\bf Proof of Theorem 1.}  Let $\cw = \cp_d$.  For $P\in \cp_d$ put $\|P\|=\sup_X|P|$.
Note that $\|P\|=0$ only when $P=0$ by the hypothesis that $X$ is non-algebraic.

A basis for $\cw$ is the set of monomials $\z^\a =\z_1^{\a_1}\z_2^{\a_2}\cdots\z_n^{\a_n}$ with
$|\a|=\a_1+\cdots +\a_n \leq d$.  Let $\cv$ be the subspace of $\cw$ spanned by 
$\{\z^\a : |\a|>0\}$.    Let $x$ denote the constant polynomial 1.  Define $L$ to be the linear
functional $L:\cw\to \bbc$ given by $L(P)=P(0)$.
Then $\l_d=\|L\|$ in the dual space of $\cw$.  Put $\d = $ the distance from 1 to $\cv$, so
$$
\d\ =\ \inf \biggl \|1-\sum_{0<|\a|\leq d}c_\a \z^\a\biggr\|.
$$
By Lemma 1, $\l_d \geq M^d$ if and only if $\|L\|\geq M^d$ if and only if ${1\over \d}\geq M^d$
if and only if $\d\leq {1\over M^d}$ if and only if 
$$
\left|  1-\sum_{k=1}^n \z_k A_k\right|\ \leq \ {1\over M^d} \quad{\rm on\ \ }X
$$
for polynomials $A_1,...,A_n\in\cp_{d-1}$.  We are done.\qed

\Theorem {2}  {\sl Let $\Sigma$ be a closed complex analytic subvariety of   $\cn$.
Let $K$ be a compact set contained in $\Sigma$.  Then $\wh K\ss\Sigma$.}

 \Note{1} In the case that $\Sigma\ss \bbc^2$ is the graph of an entire function $f$ on $\bbc$,
 with $f$ not a polynomial,  the result that $K\ss \Sigma$ implies $\wh K\ss \Sigma$ is 
 proved in Theorem 9.2 in [1].
 
 \medskip
\noindent
{\bf Proof of Theorem 2.}  Fix $a=(a_1,...,a_n)\in \cn-\Sigma$.  By definition $\l_d(a)
=\sup|P(a)|$ taken over all $P\in\cp_d$ with $\|P\|_K\leq 1$.  We make the following 
\medskip
\noindent
{\bf Claim:}
Fix $M\geq 1$.  Then $\l_d(a) >M^d$ for all large $d$.
\medskip
\noindent
{\bf Proof of Claim:}

\smallskip
\noindent
{\bf Case 1: $a=0$ so $\l_d(a)=\l_d$.}  Then $0\in \cn-\Sigma$. 
It is known that we may choose a function $H_1$, holomorphic on $\cn$ 
with $H_1=0$ on $\Sigma$ and $H_1(0)=1$.  
(See, for example, [2 Ch. VIII.A, Thm. 18].)
Put $H=1-H_1$.
Then $H$ is holomorphic on $\cn$ with $H=1$ on $\Sigma$ and $H(0)=0$.  Now $H$ has
a power series expansion
$$
H(\z)\ =\ \sum_{\a} c_\a \z^\a
$$
convergent on $\cn$.  Fix $R>0$. There exists a constant $C_R$ such that 
$$
|c_\a|\ \leq \ {C_R \over R^{|\a|} } \fa \a
$$
Fix $d$.  Put $P_d(\z) = \sum_{|\a|\leq d} c_\a\z^\a$ and put 
$\ce_d(\z) = \sum_{|\a|>d} c_\a\z^\a$.  Note that $P_d(0)=0$.
Then $H=P_d+\ce_d$ on $\cn$.  In particular, $1=P_d+\ce_d$ on $K$, so
$$
\left|1-P_d\right|\ \leq \ \left|\ce_d\right| \quad{\rm on\ \ } K.
\eqno{(1)}
$$

Fix $R_0$ so that $K\ss \{\z : |\z_k|\leq R_0 \ {\rm for\ } 1\leq k\leq n\}$. Choose $\z\in K$.  Then
$$
|\ce_d(\z)|\ \leq \sum_{|\a|>d} |c_\a||\z^\a|\ \leq \  \sum_{|\a|>d} |c_\a| R_0^{|\a|}.
$$
Fix $R$  such that ${R_0\over R}\leq {1\over 2M}$. Since $|c_\a|\leq {C_R\over R^{|\a|}}$,  we
have
$$
\sum_{|\a|>d} |c_\a| R_0^{|\a|}\ \leq \ C_R \sum_{|\a|>d} \left({R_0\over R}\right)^{|\a|}
$$
Put $t={R_0\over R}$. Then
$$
\sum_{|\a|>d} t^{|\a|}\ =\ \sum_{k=d+1}^\infty \sum_{|\a|=k} t^{|\a|}\ =\ 
\sum_{k=d+1}^\infty \left(\matrix { n+ k - 1 \cr n-1 \cr }\right) t^k\ \leq \ {1\over M^d} \quad{\rm for\ large\ \ }d.
$$
Hence,
$$
\left|\ce_d(\z)\right|\ \leq \ {C_R\over M^d}\ \leq \ {1\over (M/2)^d} \quad{\rm for\ large\ \ }d.
\eqno{(2)}
$$
From (1) and Theorem 1 we now get 
$$
\l_d\ \geq \ \left({M\over2}\right)^d  \quad{\rm for\ large\ \ }d.
$$
Since $M$ was arbitrary, the claim is proved Case 1.

\medskip
\noindent
{\bf Case 2: $a \in \cn-\Sigma$.}
Let $\chi:\cn\to \cn$ be the translation  $\chi(\z)=\z-a$.  Then $\chi(a)=0$ and $\chi$ takes the space 
$\cp_d$ isomorphically onto itself.  Put $K' =\chi(K)$ and $\Sigma' = \chi(\Sigma)$.
Since $a\notin \Sigma$,  $0\notin \Sigma'$ and so by Case 1 there exists $P'\in\cp_d$ with
$\|P'\|_{K'}\leq 1$, $|1-P'|\leq{1\over M^d}$ on $K'$ and $P'(0)=0$.  

By Case 1, $\l_d\geq M^d$ for large $d$, so there exists a polynomial $Q'$ in $\cp_d$
with $\|Q'\|_{K'}\leq 1$ and $|Q'(0)|\geq M^d$.
Now put $Q=Q'\circ \chi$.  Then $Q\in \cp_d$, $|Q(a)|\geq M^d$ and $\|Q\|_K\leq 1$.
so $\l_d(a) \geq M^d$.
Case 2 is done and so the claim holds.

 If $a\in\cn-\Sigma$, then $a$ does not belong to $\wh K$. So $\wh K\ss \Sigma$.
We are done.\qed

\Cor 1 {\sl Let $V$ be a closed complex submanifold of dimension one in $\cn$ and $K\ss V$ a 
compact subset.  Then $\wh K$ is a closed subset of $\bbc^n$.}

\Note {}  It is not true in general that $K$ compact in $\bbc^n$ implies that $\wh K$ is closed in 
$\cn$.

\medskip
\noindent
{\bf Proof of Corollary 1.}  In the notation of [1], if $K$ is a compact set in $\cn$, $\cs_K$ is 
the family of functions $\vf={1\over d}\log |P|$, $P\in \cp_d$, $|P|\leq 1$ on $K$. We define 
the extremal function
$$
\L_K(x) \ =\ \sup_{\vf\in\cs_K} \vf(x).
$$
Then for $x\in \cn$, $\L_K(x) <\infty \iff  x\in \wh K$.  That is, $\wh K \cap \cn = \{x\in\cn : \L_K(x)<\infty\}$.
Since $\L_K\equiv \infty$ in $\cn-V$, we conclude from Theorem 7.3 in [1] that in every connected
component of $V_K$ either $\L_K\equiv \infty$ or $\L_K$ is a locally bounded harmonic function.

Suppose now that $\{x_\nu\}$ is a sequence of points in $\wh K$ converging to a point $x$.
Then $x\in V$.   We must show that $x\in \wh K$.

We may assume that $x\notin K$, so $x$ lies in some connected component $\O$ of $V-K$
For large $\nu$, $x_\nu\in \O$. Since $\L_K(x_\nu)<\infty$, we have by Theorem 7.3 of [1] quoted
above, that $\L_K$ is locally bounded on $\O$.  In particular, $L_K(x)$ is finite, and so $x\in \wh K$.\qed

\Note{2} The projective hull, as introduced in [1], is a subset $\wh K\ss\bbp^n$ associated to any compact set $K\ss\bbp^n$.  In an affine chart $\cn = \bbp^n-\bbp^{n-1}$, the set  
$\wh K \cap\cn$   agrees with the definition given above. In what follows $\wh K$ will refer to
this full projective hull.

\Cor 2  {\sl  Let $D\ss \bbp^n$ be a complex hypersurface in complex projective space, and
suppose $\Sigma\ss \bbp^n-D$ is a closed subvariety of the complement.  If $K$ is 
a compact subset of $\Sigma$, then $\wh K \cap (\bbp^n-D)\ss \Sigma$. Moreover, if  $\Sigma$ is smooth of dimension
one, then $\wh K  \cap (\bbp^n-D)$ is closed in $\bbp^n-D$.}

\pf The divisor $D$ can be realized as a hyperplane section $D=v(\bbp^n)\cap \bbp^{N-1}$
under a Veronese embedding $v:\bbp^n\to \bbp^N$.  By Proposition 3.2 in [1] we know that
$\wh{v(K)} = v(\wh K)$. We now apply Theorem 2 to $v(K)\ss v(\Sigma)\ss \bbc^N=\bbp^N-\bbp^{N-1}$,
and the first statement follows.
If $\Sigma$ is a smooth curve, then   $\wh{v(K)}$ is closed in $\bbc^N$ by Corollary 1.
Hence,     $v^{-1}(\wh{v(K)})=v^{-1}(v(\wh{ K}))=\wh{K}$ is closed in $\bbp^n-D$.
\qed

\Cor 3  {\sl  Let $K\ss \Sigma\ss \bbp^n-D$ be as in Corollary 2. Suppose 
$\wh K \cap (\bbp^n-D) \ss\ss \bbp^n-D$.  Then $\wh K\cap  D=\emptyset$.
}
\pf
Suppose there exists a point $x\in \wh K\cap D$.  Then by Theorem  11.1  in [1], there exists 
a positive current $T$ of bidimension (1,1) with 
$$
\supp T \ \ss\ \wh K^- \and dd^c T \ =\  \mu - \d_x
\eqno{(3)}
$$
where $\mu $ is a probability measure on $K$ and $\wh K^-$ denotes the closure of $K$.
Let $U,V\ss\bbp^n$ be disjoint open subsets with $\wh K \cap (\bbp^n-D) \ss U$ and 
$D\ss V$.  Then by our hypothesis and the first part of (3), we have $T= \chi_U T + \chi_V T$. 
We conclude from the second part of (3) that $dd^c( \chi_U T )= \mu$ and $dd^c( \chi_V T)=-\d_x$,
both of which are impossible since $(dS, 1)=0$ for any current $S$ of dimension $2n-1$.\qed

\Cor{4} {\sl  Suppose $K=\g$  is a disjoint union of real closed curves contained in
a smooth complex analytic (non-algebraic) curve $\Sigma \ss\bbp^n-D$  with $D$ as above.
 Suppose that  $\Sigma-\g$ has only one unbounded component.
Then $\wh \g-\g$ is exactly the union of the bounded components of $\Sigma-\g$.}

\pf 
By Corollary 2, $\wh \g\cap(\bbp^n-D) \ss \Sigma$. As seen in the proof of Corollary 1, 
if one point   $x\in \Sigma-\g$ belongs to $\wh \g$, then the entire connected component of
$x$ in $\Sigma-\g$ belongs to $\wh \g$. By [1, Cor. 3.4] every bounded component
of $\Sigma-\g$ lies in $\wh\g$.  By Sadullaev's Theorem (cf. [1, \S 7]),  if the unbounded
component of $\Sigma-\g$ lies in $\wh\g$, then $\Sigma$ is algebraic.  Thus $\wh\g\cap (\bbp^n-D)$
is exactly the union of the bounded components of $\Sigma-\g$. We now apply Corollary 3.\qed

\vfill\eject
\Note{3}  We can find applications of Corollary 4 by considering graphs of 
certain functions in the disk algebra.

Let $\G$ be the unit circle in the $\z$-plane and let $A(\G)  = \{\vf\bigr|_\G : \vf $ is analytic 
in $|\z|<1$ and continuous in $|\z|\leq1\}$.  Let $K$ denote the graph of $\vf$ over $\G$, so
$$
K\ =\ \{(\z,\vf(\z)) : \z\in \G\}\ \ss\ \bbc^2.
$$
In  the following examples the   projective hull of $K$ in $\bbp^2$ is
$\wh K=\{(\z,\vf(\z)) : |\z|\leq1\}$.
 
 \smallskip 
 \item{(1)}  Let $\vf(\z) = \log(\z-2)$.  Then $K=\{(\z, \log(\z-2)):\z\in\G\}$, so $K$ is contained
 in the closed complex submanifold $\Sigma = \{(\z, w)  : e^w=\z-2\}$  of $\bbc^2$.

 \smallskip 
 \item{(2)}  Let $\vf(\z) = \sqrt{e^{2\z} + 5\z +20}$ and  $\Sigma = \{(\z,w)  : w^2 = e^{2\z}+5\z+20\}$.
 The graph $K$ of $\vf$ over $\G$ is contained in $\Sigma$ 
 
 \smallskip 
 \item{(3)}  Let $\vf(\z)$ be the restriction to $\G$ of a meromorphic function $F$ defined on all
 of $\bbc$ whose poles lies outside the unit disk. Write $F=E_1/E_2$ 
 where $E_1$ and $E_2$ are entire functions without common zeros. Then $\Sigma=\{(\z,  w): E_2(\z) w=E_1(\z)\}$ (the graph of $F$) is a closed submanifold of $\bbc^2$ which contains $K$

 \smallskip 
 \item{(4)}  Let 
 $$
 \vf(\z)  \ =\  \sum_{\nu=1}^\infty {c_\nu \over \z-a_\nu}
 $$   
 where $\sum_\nu |c_\nu| <\infty$,  each $|a_\nu|>1$,  and $\lim_{\nu\to\infty}a_\nu = a$ with 
 $|a|>1$. Then $K$ is contained in $\Sigma = \{(\z,\vf(\z)):\z\neq a {\rm \ or\ }a_\nu$ for any $\nu\}$, which is a closed submanifold
 of $\bbc^2-\{\z=a\} = \bbp^2-\bbp^1_{\infty}-\{\z=a\}$.

\vskip.3in
\centerline{\bf 
Appendix: Proof of the Banach Distance Formula.}\bigskip

Fix  $f\in\cv$.  Then $L(x-f) = L(x)=1$.  Hence, $1=|L(x-f)|\leq \|L\|\cdot\|x-f\|$.
Hence,  $1\leq \|L\|\cdot \delta$, so $\|L\|\geq {1\over \d}$.

Next, fix $g\in\cw$, $g=f+tx$ for $t\in \bbc-\{0\}$, $f\in \cv$. 
 We have $L(g)=t$.  Now $g=t(x+{f\over t})$, so 
 $\|g\| = |t|\, \|x+{f\over t} \| \geq |t|\cdot \d$, since $-{f\over t} \in \cv$. So 
 $|L(g)| = |t| \leq {1\over\d} \|g\|$.   Thus, $\|L\|\leq {1\over \d}$, and therefore 
$\|L\|  =  {1\over \d}$.\qed

\vskip.3in

\centerline{\bf 
References.}\bigskip

\item{[1]} F. R. Harvey and H. B. Lawson, Jr. {\sl Projective hulls and the projective Gelfand transform},
Asian J. Math. {\bf 10} (2006), 607-646.
\smallskip

\item{[2]}   R. C. Gunning and H. Rossi, { Analytic Functions of Several Complex Variables},
Prentice-Hall, Englewood Cliffs, N. J., 1965.

 \end

{{\parindent= .5in
\narrower\abstractfont \noindent
We prove that if $K$ is a compact subset of an affine variety $\Omega = \bbp^n-D$ ($D$ a projective hypersurface) and if  $K\ss V$, a closed analytic subvariety of $\Omega$, then the projective
hull $\wh K$ satisfies $\wh K\cap \O\ss V$.
If $V$ is smooth and 1-dimensional, then $\wh K\cap\O$ is closed in $\O$.  The result has 
applications to graphs in $\bbc^2$ of functions in the disk algebra.

}}
\vskip .4in

\end